%
%
%
%
%

\magnification=\magstephalf
\input amstex
\documentstyle{amsppt}
\nologo

\def\re{\operatorname{Re}}
\def\im{\operatorname{Im}}
\def\eps{\varepsilon}
\def\tild{\widetilde}
\def\CC{\Bbb C}
\def\RR{\Bbb R}
\def\CP{{$\Cal P$}}
\def\QED{\hfill {$\square$}}
\def\skipaline{\vskip12pt}

\topmatter

\title
On extremal mappings in complex ellipsoids
\endtitle

\author
Armen Edigarian
\endauthor

\address
\noindent
Instytut Matematyki\newline
Uniwersytet Jagiello\'nski\newline
Reymonta 4\newline
30-059 Krak\'ow, Poland
\endaddress

\email
edigaria\@im.uj.edu.pl
\endemail

\abstract
Using a generalization of \cite{Pol} we present a description of complex
geodesics in arbitrary complex ellipsoids.
\endabstract

\endtopmatter

\document

\head
1. Introduction and the main results
\endhead

Let $\Cal E(p):=\{|z_1|^{2p_1}+\dots+|z_n|^{2p_n}<1\}\subset\CC^n$,
where $p=(p_1,\dots,p_n)$, $p_j>0$, $j=1,\dots,n$; $\Cal E(p)$ is called a
{\it complex ellipsoid}.

The aim of the paper is to characterize complex $\varkappa_{\Cal E(p)}$-
and  $\tild k_{\Cal E(p)}$-geodesics. The case where $\Cal E(p)$ is convex
(i\.e\. $p_1,\dots,p_n\geq1/2$) has been solved in \cite{Jar-Pfl-Zei}.
The paper is inspired by methods of \cite{Pol}.

\skipaline

Let $D\subset\CC^n$ be a domain and let $\varphi\in\Cal O(E,D)$, where $E$
denotes the unit disk in $\CC$ and $\Cal O(\Omega,D)$ is the set of all
holomorphic mappings $\Omega\longrightarrow D$. Recall that $\varphi$ is
said to be a {\it $\varkappa_D$-geodesic \/} if there exists
$(z,X)\in D\times\CC^n$ such that:

$\varphi(0)=z$ and $\varphi'(0)=\lambda_\varphi X$ for some
$\lambda_\varphi>0$,

for any $\psi\in\Cal O(E,D)$ such that $\psi(0)=z$ and
$\psi'(0)=\lambda_\psi X$ with $\lambda_\psi>0$, we have
$\lambda_\psi\leq\lambda_\varphi$;

\vskip2pt

We say that $\varphi$ is a {\it  $\tild k_D$-geodesic \/} if there exists
$(z,w)\in D\times D$ such that:

$\varphi(0)=z$ and $\varphi(\sigma_\varphi)=w$ for some
$\sigma_\varphi\in(0,1)$,

for any $\psi\in\Cal O(E,D)$ such that $\psi(0)=z$ and
$\psi(\sigma_\psi)=w$ with $\sigma_\psi>0$,
we have $\sigma_\varphi\leq\sigma_\psi$; cf\. \cite{Pan}.

\skipaline

Let us fix some further notations:

$H^\infty(\Omega,\CC^n):=$ the space of all bounded holomorphic mappings
$\Omega\longrightarrow\CC^n$;

$\|f\|_\infty:=\sup\{\|f(z)\|\: z\in\Omega\}$, $f\in H^\infty(\Omega,\CC^n)$,
where $\|\;\|$ denotes the Euclidean norm in $\CC^n$;

$f^\ast(\zeta):=$ the non-tangential boundary value of $f$ at
$\zeta\in\partial E$, $f\in H^\infty(E,\CC^n)$;

$\Cal A(\Omega,\CC^n):=\Cal C(\bar\Omega,\CC^n)\cap\Cal O(\Omega,\CC^n)$;

$z\bullet w:=\sum_{j=1}^n z_jw_j,\quad
z\cdot w :=(z_1w_1,\dots, z_nw_n),\quad
z=(z_1,\dots,z_n),\; w=(w_1,\dots,w_n)\in\CC^n$;

$A_\nu:=\{z\in\CC\: \nu<|z|<1\}$, $\nu\in(0,1)$;

$PSH(\Omega):=$ the set of all plurisubharmonic functions on $\Omega$.

\skipaline

Fix $w_1,\dots,w_N\in\Cal A(A_\nu,\CC^n)$ and define
$$
\Phi_j(h)=\frac1{2\pi}\int_0^{2\pi}\re\Big(h^\ast(e^{i\theta})\bullet
w_j(e^{i\theta})\Big)d\theta,\quad h\in H^\infty(E,\CC^n),\quad j=1,\dots,N.
$$
We say that the functionals $\Phi_1,\dots,\Phi_N$ are {\it linearly
independent\/} if for arbitrary $s=(s_1,\dots,s_n)$, $g\in H^\infty(E,\CC^n)$,
and $\lambda_1,\dots,\lambda_N\in\RR$ such that $s_k$ nowhere vanishes on $E$,
$k=1,\dots,N$, and $g(0)=0$ the following implication is true:

if $\sum_{k=1}^N\lambda_kw_k\cdot s^\ast=g^\ast$ on a subset of $\partial E$
of positive measure, then $\lambda_1=\dots=\lambda_N=0$.

{\it Later on, we always assume that the functionals $\Phi_1,\dots,\Phi_N$
are linearly independent.}

\skipaline

{\bf Problem} (\CP).
Given a bounded domain $D\subset\CC^n$ and numbers $a_1,\dots,a_N\in\RR$,
find a mapping $f\in\Cal O(E,D)$ such that $\Phi_j(f)=a_j$, $j=1,\dots,N$,
and there is no mapping $g\in\Cal O(E,D)$ with

$\Phi_j(g)=a_j$, $j=1,\dots,N$,

$g(E)\subset\subset D$.

\skipaline

Any solution of (\CP) is called an {\it extremal mapping for (\CP)} or,
simply, an {\it  extremal}.

Problem (\CP) is a generalization of Problem (P) from \cite{Pol}.

\skipaline

We say that problem (\CP) is of type ($\Cal P_m$) if there exists a
polynomial $Q(\zeta)=\prod_{k=1}^m(\zeta-\sigma_k)$ with $\sigma_1,\dots,
\sigma_m\in E$ such that $Qw_j$ extends to a mapping of the class
$\Cal A(E,\CC^n)$, $j=1,\dots,N$.

\skipaline

One can prove that (for bounded domains $D\subset\CC^n$) any complex
$\varkappa_D$- or $\tild k_D$-geodesic may be characterized as an extremal
for a suitable problem of type ($\Cal P_1$), cf\. \S\;4.

The main result of the paper is the following

\proclaim{Theorem 1} Let $D\subset\subset G\subset\subset\CC^n$ be domains
and let $u\in PSH(G)\cap\Cal C(G)$ be such that $D=\{u<0\}$, $\partial
D=\{u=0\}$. Suppose that $f\in\Cal O(E,D)$ is an extremal for (\CP).
Assume that there exist: a set $S\subset\partial E$, a mapping $s=(s_1,
\dots,s_n)\in H^{\infty}(E,\CC^n)$, a number $\eps>0$, and a function
$v\:S\times\Cal A(E,\CC^n)\longrightarrow\CC$ such that:

{\rm (a)} $\partial E\setminus S$ has zero measure,

{\rm (b)} $f^\ast(\zeta)$, $\nabla u(f^\ast(\zeta))$, $s^\ast(\zeta)$
are defined for all $\zeta\in S$,

{\rm (c)} $s_k$ nowhere vanishes on $E$, $k=1,\dots,n$,

{\rm (d)} $u\big(f^\ast(\zeta)+s^\ast(\zeta)\cdot h(\zeta)\big)=u
\big(f^\ast(\zeta)\big)+2\re\bigg(\nabla u(f^\ast(\zeta))\bullet
\big(s^\ast(\zeta)\cdot h(\zeta)\big)\bigg)+v(\zeta,h)$,

\hfill $\zeta\in S,\quad h\in\Cal A(E,\CC^n),\;\|h\|_\infty\leq\eps$,

{\rm (e)} $\lim_{h\to0}\sup\{|v(\zeta,h)|\:\zeta\in S\}/\|h\|_\infty=0$.

Then
$$
f^\ast(\zeta)\in\partial D\quad\text{ for a.a. } \zeta\in\partial E
$$
and there exist $\varrho\in L^\infty(\partial E)$, $\varrho>0$,  $g\in
H^\infty(E,\CC^n)$, and $(\lambda_1,\dots,\lambda_N)\in\RR^N\setminus\{0\}$
such that
$$
\bigg(\sum_{k=1}^N\lambda_k w_k(\zeta)\cdot s^\ast(\zeta)\bigg)+g^\ast(\zeta)=
\varrho(\zeta) s^\ast(\zeta)\cdot\nabla u(f^\ast(\zeta))
\quad\text{ for a.a. } \zeta\in\partial E.
$$
\endproclaim

\proclaim{Remark 2}\rm Under the assumptions of Theorem 1, if
$u\in\Cal C^1(G)\cap PSH(G)$, then one can take $s:\equiv(1,\dots,1)$.
\endproclaim

\proclaim{Corollary 3} Under the assumptions of Theorem 1, if $f$ is an
extremal for $(\Cal P_m)$, then there exist $\varrho\in L^\infty(\partial E)$,
$\varrho>0$, and $g\in H^\infty(E,\CC^n)$ such that
$$
g^\ast(\zeta)=
Q(\zeta)\varrho(\zeta) s^\ast(\zeta)\cdot\nabla u(f^\ast(\zeta))
\quad\text{ for a.a. } \zeta\in\partial E.
$$
\endproclaim

Theorem 1 generalizes Theorems 2,3 in \cite{Pol} (cf\. Remark 2).
The proofs of Theorem 1 and Corollary 3 will be presented in \S\;2 and \S\;3,
respectively.

Corollary 3 give a tool to describe the extremal mappings for ($\Cal P_m$)
in the case where $D$ is an arbitrary complex ellipsoids $\Cal E(p)$.

\proclaim{Theorem 4} Let $\varphi:E\longrightarrow\Cal E(p)$ be an extremal
for ($\Cal P_m$) such that $\varphi_j\not\equiv0$, $j=1,\dots,n$. Then
$$
\varphi_j(\lambda)=
a_j\prod_{k=1}^m\Big(\frac{\lambda-\alpha_{kj}}{1-\bar\alpha_{kj}
\lambda}\Big)^{r_{kj}}
\Big(\frac{1-\bar\alpha_{kj}\lambda}{1-\bar\alpha_{k0}\lambda}
\Big)^{1\slash p_j},\quad j=1,\dots,n,
$$
where

$a_1,\dots,a_n\in \CC\setminus\{0\}$,

$\alpha_{kj}\in\bar E$, $k=1,\dots,m$, $j=0,\dots,n$,

$r_{kj}\in\{0,1\}$ and, if $r_{kj}=1$, then $\alpha_{kj}\in E$,

$\sum_{j=1}^n|a_j|^{2p_j}\prod_{k=1}^m
(\zeta-\alpha_{kj})(1-\bar\alpha_{kj}\zeta)=\prod_{k=1}^m
(\zeta-\alpha_{k0})(1-\bar\alpha_{k0}\zeta)$,\; $\zeta\in E$.

In particular, if $\varphi$ is a complex $\varkappa_{\Cal E(p)}$- or
$\tild k_{\Cal E(p)}$-geodesic, then $\varphi$ is of the above form
with $m=1$.
\endproclaim

Theorem 4 generalizes \S\;6 of \cite{Pol} and Theorem 1 in \cite{Jar-Pfl-Zei}.
The proof of Theorem 4 will be given in \S\S\;3, 4.

\proclaim{Remark 5}\rm  In the case where $\Cal E(p)$ is convex
any mapping described in Theorem~4 with $m=1$ is a complex geodesic in
$\Cal E(p)$ (\cite{Jar-Pfl-Zei}). This is not longer true if $\Cal E(p)$
is not convex, cf\. \cite{Pfl-Zwo} for the case $n=2, p_1=1, p_2<1/2$.
\endproclaim

\head{2. Proof of Theorem 1}\endhead

Note that there are two possibilities:
either $u\circ f^\ast=0$ a.e. on $\partial E$ or there
exists $\tau>0$ such that the set $\{\theta:\
u(f^\ast(e^{i\theta}))<-\tau\}$ has positive measure.
If there exists such $\tau$ lets fix one of them.
We put $$P_0:=\cases
\varnothing \quad &\text{in the first case}\\
\{\theta:\ u(f^\ast(e^{i\theta}))<-\tau\}\quad &\text{in the second case}
\endcases,$$
$A_0:=[0,2\pi)\setminus{P_0}$, and
$$p_s(h):=\frac1{2\pi}\int_{A_0}
\bigg[\re\Big(s^\ast(e^{i\theta})\cdot\nabla u(f^\ast(e^{i\theta}))\bullet
h(e^{i\theta})\Big)\bigg]^+d\theta\quad\text{ for }
h\in L^1(\partial E,\Bbb C^n),$$
where $L^1(\partial E,\Bbb C^n)$ denote the space
of all Lebesgue integrable mappings  $\partial E\to\Bbb C^n$.

\proclaim{Remark 6}\rm
(a) Under the assumptions of Theorem 1, there exists $M>0$ such that
$$\|s^{\ast}(\zeta)\cdot\nabla u(f^{\ast}(\zeta))\|\le M\quad
\text{ for a.a. }\zeta\in\partial E$$

(b) $p_s(h)$ is a seminorm on $H^1(E,\Bbb C^n)$ and
$p_s(h)\le M\|h\|_1$, where $H^1(E)$ denote the first Hardy space
of holomorphic functions,
$$H^1(E,\Bbb C^n):=\{(f_1,\dots,f_n):\ f_j\in H^1(E)\},$$
$\|\ \|_1$ denotes the norm in $H^1(E,\Bbb C^n)$.
\endproclaim

The proof of Theorem 1 is based on the following  result.
\proclaim\nofrills{Lemma 7}\;{\rm(cf. \cite{Pol}, Lemma 6).}
Under the assumptions of Theorem 1 there exist $T>0$,
$j\in\{1,\dots,N\}$, and $\delta\in\{-1,1\}$
such that 
$$\delta\Phi_j(s\cdot h)\le Tp_s(h)$$ 
for
$h\in X_j:=\{h\in H^1(E,\Bbb C^n):\ \Phi_l(s\cdot h)=0,\
l\not=j\}$.
\endproclaim

Let us for a while assume that we have Lemma 8.

\demo{Proof of Theorem 1}
By Lemma 7 it follows that
there exist $T>0$, $\delta\in\{-1,1\}$, and $j\in\{1,\dots,N\}$, such
that
$$ \delta \Phi_j(sh)\le Tp_s(h)\quad\text{ for }h\in X_j.$$
Let $\tild\Phi(h):=\delta\Phi_j(s\cdot h)$, $h\in X_j$.
Using the Hahn-Banach theorem we can extend
$\tild\Phi$ on
$L^1(\partial E,\Bbb C^n)$ (we denote this extension by $\Phi$),
in such way, that
$$\Phi(h)\le Tp_s(h)\qquad\text{for}\ h\in L^1(\partial E,\Bbb C^n).$$
We know that $p_s(h)\le M|h|_1$, where $|h|_1$ denotes 
the norm in $L^1(E,\Bbb C^n)$.
So $\Phi$ is continuous on $L^1(\partial E,\Bbb C^n)$.
By Riesz's theorem, $\Phi$ can be represented as
$\Phi(h)=\frac1{2\pi}
\int_0^{2\pi}\re(h^\ast(e^{i\theta})\bullet\tild w(e^{i\theta}))d\theta$,
where $\tild w\in L^\infty(\partial E,\Bbb C^n)$.\par
It is easy to see that there are
$\lambda_1,\dots,\lambda_N$, not equal
simultaneously zero, such that
$\Phi(h)=\sum_{k=1}^N\lambda_k\Phi_k(s\cdot h)$ for
$h\in H^1(E,\Bbb C^n)$.
We denote by $G$ the linear functional on $L^1(\partial E,\Bbb C^n)$ defined
by the formula
$$ G(h):=\frac1{2\pi}
\int_0^{2\pi}\re\Big(\sum_{k=1}^N\lambda_k
w_k(e^{i\theta})\bullet s^\ast(e^{i\theta})\cdot
h(e^{i\theta})\Big)d\theta.$$
Then $\Phi(h)-G(h)=0$ for $h\in H^1(E,\Bbb C^n)$. By the theorem of
F\. \& M\. Riesz it
follows that there exists $g\in H^{\infty}(E,\Bbb C^n)$,  $g(0)=0$,
such that
$$\tild w-s^\ast\cdot\sum_{k=1}^N\lambda_k w_{k}=g^\ast.$$
We have
$$\multline
\Phi(h)=\frac1{2\pi}\int_0^{2\pi}\re
\Big[\Big(\sum_{k=1}^N\lambda_kw_k(e^{i\theta})\cdot s^\ast(e^{i\theta})
+g^\ast(e^{i\theta})\Big)\bullet h^\ast(e^{i\theta})\Big]d\theta\\
\le T\frac1{2\pi}\int_{A_0}\Big[\re\Big(\big(s^\ast(e^{i\theta})\cdot
\nabla u(f^\ast(e^{i\theta})\big)\bullet h^\ast(e^{i\theta})\Big)\Big]^+
d\theta
\endmultline
\tag 1
$$
for any $h\in H^1(E,\Bbb C^n)$. We see that right-hand side is zero
for any $h\in H^1(E,\Bbb C^n)$ (hence, for any $h\in L^1(E,\Bbb C^n)$)
such that
$$h^\ast\le0\quad\text{ on }\partial E\setminus{\Big( P_0\cup
\{\zeta\in\partial E:\ s^\ast(\zeta)\cdot\nabla u(f^\ast(\zeta))=0\}\Big)}.$$

Hence
$$\sum_{k=1}^N\lambda_kw_k\cdot s^\ast+g^\ast=0\quad\text{ a.a. on }
P_0\cup\{\zeta\in\partial E:\ s^\ast(\zeta)\cdot\nabla u(f^\ast(\zeta))=0\}.$$
We know that $\Phi_1,\dots,\Phi_N$ are linearly independent, so
Lebesgue measure of $P_0$ and
$\{\zeta\in\partial E:\ s^\ast(\zeta)\cdot\nabla u(f^\ast(\zeta))=0\}$ are
equal zero.
Hence
$$\sum_{k=1}^N\lambda_kw_k\cdot s^\ast+g^\ast=
\varrho s^\ast(\zeta)\cdot\nabla u(f^\ast(\zeta)),$$ where
$\varrho(\zeta)\in\Bbb C\setminus{\{0\}}$ for a.a. $\zeta\in\partial E$.
Now, it is enough to remark that condition (1) implies that
$0<\varrho\le T$ a.a. on $\partial E$.
\QED
\enddemo
Now, we are going to prove Lemma 7.
\demo{Proof of Lemma 7}
Suppose that the lemma is not true. Then for each
 $j\in\{1,\dots,N\}$ and $m\in\Bbb N$ there are
$h_{jm}^{+},h_{jm}^{-}\in X_j$,
such that
$$\Phi_j(s\cdot h_{jm}^{+})> mp_s(h_{jm}^{+}), \quad
-\Phi_j(s\cdot h_{jm}^{-})> mp_s(h_{jm}^{-}).$$
We may assume that $h^+_{jm}, h^-_{jm}\in\Cal A(E,\Bbb C^n)$ and that
$$\Phi_j(s\cdot h_{jm}^{+})=1,\quad \Phi_j(s\cdot h_{jm}^{-})=-1.$$
For any
$q=(q_1^+,q_1^-,\dots,q_N^+,q_N^-)\in\Bbb R^{2N}_+$ we define
the function
$$f_{qm}=f+\sum_{j=1}^N(q_j^+s\cdot h_{jm}^++q_j^-s\cdot h_{jm}^-)=
f+s\cdot h_{qm}$$
and the linear mapping of $A:\Bbb R_+^{2N}\to\Bbb R^{N}$,
$A(q):=(q_1^+-q_1^-,\dots,q_N^+-q_N^-)$.
Note that $\Phi_j(f_{qm})-\Phi_j(f)=A(q)_j$.


\proclaim\nofrills{Lemma 8}\;{\rm(see \cite{Pol}, Lemma 7).}  
Let $u$ be a non-positive subharmonic
function in $E$ and let $\triangle u$ be the Riesz measure of $u$.
Suppose that one of the following conditions is true
\roster
\item"(a)"
$\triangle u(r_0E)>a>0$ for some $r_0\in(0,1)$,
\item"(b)" for some set
$Z\subset[0,2\pi)$ with positive measure, the upper radial limits
of $u$ at $\zeta\in Z$ do not exceed $-a<0$
\rom(i.e. $\limsup_{r\to 1}u(r\zeta)\le -a$\rom).
\endroster
Then $u(\zeta)\le-C(1-|\zeta|)$, where $C>0$ is a constant
depending only on $r_0$, $a$, and $Z$.
\endproclaim

\skipaline

Let $u_0:=u\circ f$.
\proclaim{Lemma 9} There exist constant $C>0$ and
constants $t_m>0$, $m\in\Bbb N$, such that for $\|q\|<t_m$ we have
\roster
\item"(a)" $f_{qm}\in\Cal O(E,G)$ (so, we define $u_{qm}:=u\circ f_{qm}$),
\item"(b)" $u_{qm}(\zeta)\le v_{qm}(\zeta):=C\ln|\zeta|+
\frac1{2\pi}\int_{A_0}\big[u^\ast_{qm}(e^{i\theta})\big]^+
P(\zeta,\theta)d\theta$ for $|\zeta|>\frac12$.
\endroster
\endproclaim
\demo{Proof of Lemma 9}
Ad (a). It follows from the assumption that $D\subset\subset G$.\par
Ad (b). Suppose  that there exists $r_0\in(0,1)$  such that
$\triangle u_0(r_0E)>a>0$.  The continuity of $u$ implies that for
$$\tild u_{qm}(\zeta):=u_{qm}(\zeta)-
\frac1{2\pi}\int_{A_0}\big[u_{qm}^\ast(e^{i\theta})\big]^+
P(\zeta,\theta)d\theta,\quad\zeta\in E,$$
if $t_m$ is small enough then $\triangle\tild u_{qm}(rE)>\frac a2$.
Hence, from Lemma 8 we get the required result.\par
If  $\triangle u_0(rE)=0$ for any $r\in(0,1)$ and
$u_0^\ast(\zeta)=0$ for a.a. $\zeta\in\partial E$, then from (2)
we get that $u_0$ is harmonic in $E$. But, it is a contradiction, since
$u_0\not\equiv0$. Hence, $P_0$ has positive measure.
From the continuity of $u$ we conclude that if $t_m$ are small enough, then
$\{\zeta:\ \tild u_{qm}(\zeta)<-\frac\tau2\}$ has positive measure.
By Lemma 8 we get the required result.
\QED
\enddemo
Let us introduce some new notation:
$E_{qm}:=\{\zeta\in E:v_{qm}(\zeta)<0\}$ and
$$g_{qm}(\zeta):=\zeta\exp\Big\{\frac1{2\pi C}
\int_{A_0}\Big[u^\ast_{qm}(e^{i\theta})\Big]^+
S(\zeta,\theta) d\theta\Big\}.$$
Here $S(\zeta,\theta):=\frac{\zeta+e^{i\theta}}{\zeta-e^{i\theta}}$
is the Schwarz kernel.
\proclaim{Remark 10}\rm Note that $C\ln|g_{qm}|=v_{gm}$,
$v_{qm}(\zeta)\ge C\ln|\zeta|$ (hence, $|g_{qm}(\zeta)|\ge|\zeta|$),
and $E_{qm}=g_{qm}^{-1}(E)$.
\endproclaim
\proclaim\nofrills{Lemma 11}\;{\rm(cf\. \cite{Pol}, statement 2).}
{\rm (a)} $E_{qm}$ is connected, $0\in E_{qm}$, {\rm (b)}
$g_{qm}$ maps $E_{qm}$ conformally onto $E$.
\endproclaim
\demo{Proof of Lemma 12} (a) Note that
$E_{qm}=\bigcup_{\delta>0}\{\zeta:v_{qm}(\zeta)<-\delta\}$ and
$$\{\zeta:\ v_{qm}(\zeta)<-\delta\}\subset
 \{\zeta:\ |\zeta|<e^{-\delta\slash C}\}.$$
Since $v_{qm}$ is harmonic outside  0 and
$v_{qm}^\ast(e^{i\theta})\ge0$,
any connected component of $\{\zeta:v_{qm}(\zeta)<-\delta\}$ must contain 0.
\par
(b) At first let us see that $g_{qm}:E_{qm}\to E$ is proper.
Let $\zeta_k\to\zeta_0\in\partial E_{gm}$.
If $\zeta_0\in\partial E$, then $|g_{qm}(\zeta_k)|\to1$ (since
$|g_{qm}|\ge|\zeta|$).
If $\zeta_0\in E$,  then $|g_{qm}(\zeta_k)|\to|g_{qm}(\zeta_0)|=1$.

Since $g_{qm}'(0)\ne0$ and $g_{qm}^{-1}(0)=\{0\}$,
$g_{qm}$ is conformal.
\QED
\enddemo
We define $\tild f_{qm}(\zeta)=f_{qm}(g_{qm}^{-1}(\zeta))$,
$\widehat f_{qm}(\zeta)=\tild f_{qm}(e^{-\|q\|\slash m}\zeta)$,
$$\tild A_m(q)=(\Phi_1(\tild f_{qm})-\Phi_1(f),\dots,
\Phi_N(\tild f_{qm})-\Phi_N(f)),$$
and
$$\widehat A_m(q)=(\Phi_1(\widehat f_{qm})-\Phi_1(f),\dots,
\Phi_N(\widehat f_{qm})-\Phi_N(f)).$$
\proclaim{Remark 12}\rm It is easy to see that
$\tild f_{qm}(E)\subset D$, $\widehat f_{qm}(E)\subset\subset D$,
 and $\tild A_m(0)=\widehat A_m(0)=0$.
\endproclaim
Following result explains why we have used the functionals
of the special form.
\proclaim{Lemma 13}
Suppose that
$$\Phi(h)=\frac1{2\pi}\int_0^{2\pi}\re\Big(h^\ast(e^{i\theta})\bullet
 w(e^{i\theta})\Big)d\theta,$$
where
$w\in\Cal A(A_\nu,\Bbb C^n)$ for some $\nu\in(0,1)$,
$f\in H^\infty(E,\Bbb C^n)$, and that $g\in\Cal O(E,E)$, $g(0)=0$.
Then
$$ |\Phi(f\circ g)-\Phi(f)|\le K\|f\|_\infty\sup_{\zeta\in E}
|g(\nu\zeta)-\nu\zeta|,$$
where $K>0$ depends only on $\Phi$.
\endproclaim
\demo{Proof of Lemma 13}
We have
$$\Phi(h)=\frac1{2\pi}\int_0^{2\pi}\re\Big(h(\nu e^{i\theta})\bullet
 w(\nu e^{i\theta})\Big)d\theta.\tag 4$$
Hence
$$|\Phi(h)|\le\Big(\max_{\zeta\in\partial E}\|w(\nu\zeta)\|\Big)
\Big(\max_{\zeta\in\partial E}\|h(\nu\zeta)\|\Big).$$
But,
$$\|f(g(\nu\zeta))-f(\nu\zeta)\|\le\Big(\sup_{\xi\in E}|f'(\nu\xi)|\Big)
|g(\nu\zeta)-\nu\zeta|,$$
and $\sup_{\xi\in E}|f'(\nu\xi)|\le\frac{\|f\|_\infty}{1-\nu^2}$.
\QED
\enddemo
\proclaim\nofrills{Lemma 14}\;{\rm(cf\. \cite{Pol}, Statement 3).}
The mappings $\tild A_m, \widehat A_m$ are
continuous in $q$ when $\|q\|<t_m$.
\endproclaim
\demo{Proof of Lemma 14}
 It is enough to remark that
 if $q_k\to q$, then $u_{q_km}^\ast\to u_{qm}^\ast$
 uniformly on $\partial E$.
Hence
$g_{q_km}\to g_{qm}$ uniformly on compact sets of $E$. It is evident after
the last assertion that $g_{q_km}^{-1}\to g_{q_m}^{-1}$,
$\tild f_{q_km}\to\tild f_{qm}$, and
$\widehat f_{q_km}\to\widehat f_{qm}$
 uniformly on compact sets, too.
Since $\Phi_j$ are continuous with respect to this convergence
(it follows easily from (4)), we conclude the proof.
\QED
\enddemo
\proclaim{Lemma 15} For each $b>0$ there is $m_0\in\Bbb N$
such that for any $m\ge m_0$ there is $q_m>0$
such that $\|A(q)-\tild A_m(q)\|\le b \|q\|$, when $\|q\|\le q_m$.
\endproclaim
\demo{Proof of Lemma 15}
It follows from the definition of $A,\tild A_m$
that it is enough to prove inequality
$$|\Phi(\tild f_{qm})-\Phi(f_{qm})|\le b\|q\|,$$
where $\Phi$ is the functional  of our special form.
By Lemma 14 it is enough to consider inequality
$$\sup_{\zeta\in\nu E}|g^{-1}_{qm}(\zeta)-\zeta|\le b\|q\|.$$
Note that
$$\sup_{\zeta\in\nu E}|g^{-1}_{qm}(\zeta)-\zeta|\le
\sup_{\zeta\in\nu E}|g_{qm}(\zeta)-\zeta|$$
and for  small $q_m$ (such that $|1-\exp q_m|\le 2q_m$)
$$\multline
\Big|1-\exp\Big(\frac1{2\pi C}
\int_{A_0}\Big[u^\ast_{qm}(e^{i\theta})\Big]^+
S(\zeta,\theta) d\theta\Big)\Big|\\
\le
2\frac{1+\nu}{1-\nu}\Big(\frac1{2\pi C}
\int_{A_0}\Big[u^\ast_{qm}(e^{i\theta})\Big]^+
d\theta\Big)
\endmultline
$$
for $\zeta\in\nu E$.
Hence, it is enough to consider
$$
\multline
\int_{A_0}\Big[u_{qm}^\ast(e^{i\theta})\Big]^+d\theta\le
\int_{A_0}2\Big[\re(\nabla u(f^\ast(e^{i\theta}))\bullet
s^\ast(e^{i\theta})\cdot h_{qm}(e^{i\theta}))\Big]^+d\theta
\\ +o(\|h_{qm}\|_\infty).
\endmultline
$$
But $p_s(h_{qm})\le\|q\|\max\{p_s(h_{jm}): j=1,\dots,N\}\le\frac1m\|q\|$.
Hence, if $m$ is big and $q_m$ is small enough, we get the required result.
\QED
\enddemo

\proclaim{Lemma 16} For each $b>0$ there is $m_0\in\Bbb N$
such that for any $m\ge m_0$ there is $q_m>0$
such that $\|\tild A_m(q)-\widehat
A_m(q)\|\le b \|q\|$, when $\|q\|\le q_m$.
\endproclaim
\demo{Proof of Lemma 16} As in Lemma 15, by Lemma 13 it is enough to
proof inequality
$$\sup_{\zeta\in\nu E}|e^{-\|q\|\slash m}\zeta-\zeta|\le b\|q\|.$$
But, for a small $\|q\|\slash m$ we have
$|1-e^{-\|q\|\slash m}|\le 2\frac{\|q\|}m$. Hence, we get the required result.
\QED
\enddemo

\proclaim\nofrills{Lemma 17}\;{\rm (cf\. \cite{Pol}, Lemma 8).}
For any continuous mapping
$F:\Bbb R^{2N}_+\to \Bbb R^{N}$,
 if $$\|F(x)-A(x)\|\le b\|x\| \quad\text{ for }
  x\in  B(0,r)\cap\Bbb R_+^{2N},$$
 where $b=\frac1{2\sqrt N}$,
then there exists $q\in B(0,r)\cap R^{2N}_+\setminus{\{0\}}$ such that
$F(q)=0$.
\endproclaim
\demo{Proof of Lemma 17}
Let us denote
$$\Cal Q:=\{(x_1,\dots,x_N):\ 0<x_j<t_0,\ j=1,\dots,N\}$$
and
$$\pi:\Bbb R^{N}\ni(x_1,\dots,x_N)\to
(x_1,t_0-x_1,\dots,x_N,t_0-x_N)\in\Bbb R^{2N},$$
where
$t_0=\frac1{2\sqrt{N}}\min\{1,r\}$.
It easy to check that $\|\pi(l)\|\le t_0\sqrt N$ for $l\in\bar\Cal Q$ and
$\pi(\Cal Q)\subset B(0,r)\cap\Bbb R_+^{2N}$.
 Note, that
$$\|F\circ \pi(l)-A\circ\pi(l)\|\le
b\|\pi(l)\|\le \frac{t_0}2\quad \text{for } l\in\bar\Cal Q.$$
Let us consider the homotopy defined by the formula
$\tild F_t=t F\circ\pi+(1-t)A\circ\pi$.
It is enough to show that $0\not\in\tild F_t(\partial\Cal Q)$.
Then from the homotopical invariance of the degree of mappings \cite{Zei}
we have:  $\deg(F\circ\pi,\Cal Q,0)=\deg(A\circ\pi,\Cal Q,0)\not=0$,
hence
$0\in F\circ\pi(\Cal Q)$.\par
It is easy to see that for any $l\in\partial\Cal Q$
$$t_0\le\|A\circ\pi(l)\|\le \|\tild F_t(l)\|+t\|F\circ\pi(l)
-A\circ\pi(l)\|\le\|\tild F_t(l)\|+
\frac{t_0}2.$$
Hence, we get the required result.
\QED
\enddemo

Let us return to the proof of  Lemma 9. By Lemmas 14,15, and 16
it follows that
$\widehat A_m$ is continuous in $\Bbb R_+^{2N}$ and for each $b>0$ there is
an $m\in\Bbb N$ and $q_m>0$, such that $\|\widehat A_m(q)-A(q)\|\le b\|q\|$.
By Lemma 17, for some
$m$  we can find $q_0$, which is a solution of the
equation $\widehat A_m(q_0)=0$.
Hence, we have
$$ \Phi_j(\widehat f_{q_0 m})=a_j\quad\text{ for } j=1,\dots,N.$$
But, this contradicts the extremality of $f$, since
$\widehat f_{q_0 m}(E)\subset\subset D$.
\QED
\enddemo

\head{3. Proof of Theorem 4}\endhead

Before we prove the theorem we recall some auxiliary results.
\proclaim{Lemma 18} 
Let $\varphi\in H^1(E)$ be such that
$$\frac{\varphi^\ast(\zeta)}{\prod_{k=1}^m(\zeta-\sigma_k)}\in\Bbb R_{>0}\quad
\text{ for a.a. } \zeta\in\partial E,$$
where $\sigma_k\in\Bbb C$, $k=1,\dots,m$. Then
there exist $r\in\Bbb R$ and $\alpha_k\in\bar E$,
$k=1,\dots,m$, such that
$$\varphi(\zeta)=r\frac{\prod_{k=1}^m(\zeta-\alpha_k)
(1-\bar\alpha_k\zeta)}
{\prod_{k=1}^m(1-\bar\sigma_k\zeta)},
\quad\zeta\in E.$$
\endproclaim
This Lemma is generalization of Lemma 8.4.6 in \cite{Jar-Pfl}.
\demo{Proof of Lemma 18}
Put $\tild\varphi(\zeta)=\varphi(\zeta)\prod_{k=1}^m(1-\bar\sigma_k\zeta)$,
then $\tild\varphi\in H^1(E)$ and
$$\frac1{\zeta^m}\tild\varphi^\ast(\zeta)\in\Bbb R_{>0}\quad
\text{ for a.a. } \zeta\in\partial E.$$
Hence, it is enough to prove the
lemma for $\sigma_k=0$, $k=1,\dots,m$. Let us
denote $$P(\zeta)=\sum_{k=0}^m\frac{\varphi^{(k)}(0)}{k!}\zeta^k+
\sum_{k=0}^{m-1}\frac{{\overline{\varphi^{(k)}(0)}}}{k!}\zeta^{2m-k}.$$
It is  easy to see that if
$\psi(\zeta):=\frac{\varphi(\zeta)-P(\zeta)}{\zeta^m}$, then
$\psi\in H^1(E)$ and $\psi^\ast(\zeta)\in\Bbb R$ for a.a.
$\zeta\in\partial E$. Hence $\psi\equiv0$.

Let $t(\theta):=\frac{P(e^{i\theta})}{e^{i\theta m}}$.
We know that $t$ is $\Bbb R$-analytic, $t(\theta)\ge0$ for $\theta\in\Bbb R$.
If for some $\theta_0\in\Bbb R$ we have $t(\theta_0)=0$ then
$t(\theta)=(\theta-\theta_0)^k\tild t(\theta)$, where $k$ is even.

Note that $\overline{P(1\slash\bar\zeta)}=\frac{P(\zeta)}{\zeta^{2m}}$ and
if $P(0)=0$, then
$P(\zeta)=\zeta^k\tild P(\zeta)$, $\tild P(0)\not=0$,
$\deg\tild P=2m-2k$, and
$\overline{\tild P(1\slash\bar\zeta)}=\frac{\tild P(\zeta)}{\zeta^{2(m-k)}}$
Now, it is enough to
note that
if $P(\zeta_0)=0$, $\zeta_0\not=0$, then
$P(1\slash\bar\zeta_0)=0$ and if
$Q(\zeta):=\frac{P(\zeta)}{(\zeta-\zeta_0)(1-\bar\zeta_0\zeta)}$,
then
$\overline{Q(1\slash\bar\zeta)}=\frac{Q(\zeta)}{\zeta^{2(m-1)}}$.
\QED
\enddemo
\proclaim{Lemma 19} Let $S_1,S_2$ be singular inner
functions and let $S_1S_2\equiv1$. Then $S_1, S_2\equiv 1$.
\endproclaim
\demo{Proof of Lemma 19}
Suppose that
$S_j(z)=\exp\Big(-\int_0^{2\pi}\frac{e^{it}+z}{e^{it}-z}d\mu_j(t)\Big)$,
$j=1,2$, where $\mu_1,\ \mu_2$ are non-negative Borel measures,
singular w\.r\.t\. Lebesgue measure.
Then $S_1S_2\equiv1$ is equivalent to
$\mu_1+\mu_2=0$. Since $\mu_j\ge0$, $j=1,2$,  $\mu_1=\mu_2=0$.
\QED
\enddemo


\demo{Proof of Theorem 4}
We know that $\varphi_j=B_jS_jF_j$, where
$B_j$ is a Blaschke product, $S_j$ is a singular inner function and
$F_j$ is an outer function.
Let us take $s:=(F_1,\dots,F_n)$. Note that
$\Big|{\varphi_j^\ast(\zeta)\slash F_j^\ast(\zeta)}\Big|=1$
for a.a. $\zeta\in\partial E$
and
$\frac{\partial u}{\partial z_j}(\varphi)=
p_j\frac{|\varphi_j|^{2p_j}}{\varphi_j}$
for
$j=1,\dots,n$.
We want to show that the assumptions
of  Theorem 1 are fulfilled.
Let $u(z):=\sum_{j=1}^n|z_j|^{2p_j}-1$ be the defining function for
$\Cal E(p)$.\par
We know that $\varphi_j\not\equiv0$, $j=1,\dots,n$.
Hence $\nabla u(\varphi^\ast(\zeta))$
exists for a.a.$\zeta\in\partial E$.\par
We have
$$\multline
\frac{|\varphi_j+F_jh_j|^{2p_j}-|\varphi_j|^{2p_j}-2\re\Big(p_j\frac
{|\varphi_j|^{2p_j}}{\varphi_j}F_jh_j\Big)}{|h_j|}\\
=|\varphi_j|^{2p_j}\frac{|1+\frac{F_j}{\varphi_j}h_j|^{2p_j}-1
-2p_j\re\Big(\frac
{F_j}{\varphi_j}h_j\Big)}{|h_j\frac{F_j}{\varphi_j}|}.
\endmultline$$
From the equality
$$\lim_{z\to 0}\frac{|1+z|^\alpha-1-\alpha\re z}{|z|}=0
\qquad\alpha>0,$$
we see that all the assumptions of Theorem 1 are fulfilled.

\skipaline

Hence, by Corollary 3,
there exist $g\in H^{\infty}(E,\Bbb C^n)$ and $\varrho\in L^{\infty}(E)$,
$\varrho>0$, such that
$$Q(\zeta)\varrho(\zeta)F^\ast_j(\zeta)\frac{|\varphi^\ast_j(\zeta)|^{2p_j}}
{\varphi^\ast_j(\zeta)}
=g_j^\ast(\zeta)\quad
\text{ for a.a. } \zeta\in\partial E,\ j=1,\dots,n,$$
where $Q(\zeta)=\prod_{k=1}^m(\zeta-\sigma_k)$ is polynomial for problem
($\Cal P_m$).
It is equivalent to
$$Q(\zeta)\varrho(\zeta)|F^\ast_j(\zeta)|^{2p_j}
=B^\ast_j(\zeta)S^\ast_j(\zeta)g_j^\ast(\zeta)\quad
\text{ for a.a. } \zeta\in\partial E,\ j=1,\dots,n.$$
By Lemma 18 there exist $r_j>0$ and $\alpha_{kj}\in\bar E$ such
that
$$B^\ast_j(\zeta)S^\ast_j(\zeta)g_j^\ast(\zeta)=r_j\frac{\prod_{k=1}^m
(\zeta-\alpha_{kj})
(1-\bar\alpha_{kj}\zeta)}{\prod_{k=1}^m(1-\bar\sigma_k\zeta)}\tag 3$$
and there exist $r_0>0$ and $\alpha_{k0}\in\bar E$ such that
$$Q(\zeta)\varrho(\zeta)=\sum_{j=1}^n
B^\ast_j(\zeta)S^\ast_j(\zeta)g_j^\ast(\zeta)=
r_0\frac{\prod_{k=1}^m (\zeta-\alpha_{k0})
(1-\bar\alpha_{k0}\zeta)}{\prod_{k=1}^m(1-\bar\sigma_k\zeta)}.\tag 4$$
We have
$$r_0\prod_{k=1}^m (\zeta-\alpha_{k0})
(1-\bar\alpha_{k0}\zeta)|F_j(\zeta)|^{2p_j}=
r_j\prod_{k=1}^m (\zeta-\alpha_{kj})(1-\bar\alpha_{kj}\zeta).\tag 5$$
Hence
$$F_j(\zeta)=a_j\prod_{k=1}^m\Big(\frac{1-\bar\alpha_{kj}\zeta}
{1-\bar\alpha_{k0}\zeta}
\Big)^{1\slash p_j},\tag 6$$
where $a_j\in\Bbb C\setminus{\{0\}}$.
From (6) it follows that
$$B_j(\zeta)=\prod_{k=1}^m\Big(
\frac{\zeta-\alpha_{kj}}{1-\bar\alpha_{kj}\zeta}\Big)^{r_{kj}},
\text{ where } r_{kj}\in\{0,1\}.$$
Hence
$$S_j(\zeta)g_j(\zeta)=
r_j\frac{\prod_{k=1}^m
(\zeta-\alpha_{kj})^{1-r_{kj}}
(1-\bar\alpha_{kj}\zeta)^{1+r_{kj}}}{\prod_{k=1}^m(1-\bar\sigma_k\zeta)}.$$
Since the right hand side is an outer function,
from Lemma 19 we get that $S_j\equiv1$, $j=1,\dots,n$.\par
From (5) and (6) we see that
$|a_j|^{2p_j}=\frac{r_j}{r_0}$
 and from (3) and (4) it follows that
$\sum_{j=1}^n|a_j|^{2p_j}\prod_{k=1}^m
(\zeta-\alpha_{kj})(1-\bar\alpha_{kj}\zeta)=\prod_{k=1}^m
(\zeta-\alpha_{k0})(1-\bar\alpha_{k0}\zeta)$, $\zeta\in E$.
So, we get the required result.
\QED
\enddemo

\head{4. The case of complex geodesics}\endhead
\proclaim{Lemma 20} Any  $\varkappa_D$- and $\tild k_D$-geodesic
is  extremal for appropriate problem ($\Cal P_1$).
\endproclaim
\demo{Proof of Lemma 20}
The case of $\varkappa_D$-geodesic.
Let us consider problem (\CP) with linear functionals such that:

$N=4n$, 

$w_j:=(0,\dots,1,\dots,0)$, $a_j:=\re z_j$, for $j=1,\dots,n$, 

$w_j:=(0,\dots,-i,\dots,0)$, $a_j:=\im z_j$, for $j=n+1,\dots,2n$,

$w_j:=(0,\dots,\frac1{\zeta},\dots,0)$, $a_j:=\re X_j$, for  $j=2n+1,\dots,3n$,

$w_j:=(0,\dots,\frac{-i}{\zeta},\dots,0)$, $a_j:=\im X_j$, 
for $j=3n+1,\dots,4n$,

where $z\in D$ and $X\in\Bbb C^n\setminus{\{0\}}$.

It is easy to see that corresponding linear functionals are linearly
independent and the problem (\CP) has degree 1.\par
Let us show that any  $\varkappa_D$-geodesic for $(z,X)$
 $f$ is  extremal for this problem (\CP). Suppose that there exists a mapping
$g\in\Cal O(E,D)$ such that   $g(0)=z$, $g'(0)=X$,
and
$g(E)\subset\subset D$. Denote $\tild g(\zeta):=g(\zeta)+\zeta tX$,
where $t>0$ will be defined later. Then $\tild g(0)=g(0)=z$ and
$\tild g'(0)=g'(0)+tX=(1+t)X$. If we take $t$ such that
$\tild g(E)\subset D$ (it is possible, because $g(E)\subset\subset D$),
then we have contradiction with that $f$ is $\varkappa_D$-geodesic.\par
The case $\tild k_D$-geodesic.
 Let us consider problem (\CP) with linear functionals such that
$f\in\Cal O(E,D)$ is extremal iff
$f(0)=z$, $f(\sigma)=w$, where $\sigma>0$, and
there is no mapping $g\in\Cal O(E,D)$ such that
\roster
\item $g(0)=z$, $g(\sigma)=w$,
\item $g(E)\subset\subset D$.
\endroster
(The functions $w_j$ in this case can be
constructed by the similar way as for $\varkappa_D$-geodesic.
It is enough to replace $\frac1{\zeta}$ by
$\frac1{\zeta-\sigma}$ and
$\frac{-i}{\zeta}$ by $\frac{-i}{\zeta-\sigma}$).
It is easy to see that suitable linear functionals are linearly
independent and that we have problem ($\Cal P_1$).\par
Let us show that any  $\tild k_D$-geodesic $f$ is  extremal for
this problem (\CP). Suppose that there a exists mapping
$g\in\Cal O(E,D)$ such that   $g(0)=z$, $g(\sigma)=w$, and
$g(E)\subset\subset D$. Denote
$\tild g(\zeta):=g(\zeta)+\frac{\zeta}{t\sigma}(g(\sigma)-g(t\sigma))$,
where $0<t<1$ will be defined later. Then $\tild g(0)=g(0)=z$ and
$\tild g(t\sigma)=g(\sigma)=w$. If we take $t$ such that
$\tild g(E)\subset D$ (it is possible, because $g(E)\subset\subset D$),
then we have contradiction, because $f$ is
$\tild k_D$-geodesic.
\QED
\enddemo
\head{Acknowledgement}\endhead
The ideas of the paper comes after fruitful conversations with
Prof. M\. Jarnicki and W\. Zwonek. I would like to thank them.


\enddocument